\RequirePackage{ifpdf}
\ifpdf 
\documentclass[pdftex]{sigma}
\else
\documentclass{sigma}
\fi

\begin{document}

\renewcommand{\PaperNumber}{035}

\FirstPageHeading

\def\C{{\cal C}}
\def\L{{\cal L}}
\def\E{{\cal E}}
\def\K{{\cal K}}
\def\H{{\cal H}}
\def\U{{\cal U}}
\def\V{{\cal V}}
\def\W{{\cal W}}
\def\F{{\cal F}}
\def\M{{\cal M}}
\def\D{{\cal D}}
\def\G{{\cal G}}
\def\H{{\cal H}}
\def\A{{\cal A}}
\def\B{{\cal B}}
\def\I{{\cal I}}
\def\S{{\cal S}}

\ShortArticleName{Relative dif\/ferential $K$-characters}

\ArticleName{Relative dif\/ferential $\boldsymbol{K}$-characters}

\Author{Mohamed MAGHFOUL}

\AuthorNameForHeading{M. Maghfoul}

\Address{Universit\'e Ibn Tofa{\"{i}}l, D\'epartement de
Math\'ematiques, K\'enitra, Maroc}

\Email{\href{mailto:mmaghfoul@lycos.com}{mmaghfoul@lycos.com}}

\ArticleDates{Received November 26, 2007, in f\/inal form March
17, 2008; Published online March 28, 2008}

\Abstract{We def\/ine a group of relative dif\/ferential
$K$-characters associated with a smooth map between two  smooth
compact manifolds. We show that this group f\/its into a short
exact sequence as in the non-relative case. Some secondary
geometric invariants are expressed in this theory.}

 \Keywords{geometric $K$-homology; dif\/ferential $K$-characters}

  \Classification{51H25; 51P05; 58J28}

\section{Introduction} Cheeger and Simons \cite{CS} introduced the notion of dif\/ferential characters to express some secon\-dary geometric invariants of a principal $G$-bundle in the base space. This theory   has been appearing more and more frequently in quantum f\/ield and string theories (see \cite{BTW,LU,HS}). On the other hand, it was shown recently (see \cite{AST,Periwal,RS}) that $K$-homology of Baum--Douglas~\cite{BD1}  is an appropriate arena in which various aspects of $D$-branes in superstring theory can be described.

 In \cite{BM} we have def\/ined with M.T.~Benameur  the notion of dif\/ferential characters in $K$-theory on a smooth compact manifold. Our original motivation was to explain some secondary geometric invariants coming from the Chern--Weil and Cheeger--Simons theory in the language of  $K$-theory. To do this, we have used  the  Baum--Douglas construction of $K$-homology. As a~result,  we obtained the eta invariant of Atiyah--Patodi--Singer as a $\mathbb R/\mathbb  Z$-dif\/ferential $K$-character, while it is a $\mathbb R/\mathbb  Q$-invariant in  the works of Cheeger and Simons. Recall that a geometric $K$-cycle of Baum--Douglas  over a smooth compact manifold  $X$ is a triple $(M,E,\phi)$
 such that:
 $M$ is a~closed smooth ${\rm Spin}^{c}$-compact  manifold with a f\/ixed Riemannian structure;
 $E$ is a Hermitian vector bundle over~$M$ with a f\/ixed Hermitian connection $\nabla ^E$ and
$\phi : M \to X$ is a smooth map. Let  $\C_*(X)$ be the semi-group
for the disjoint union of equivalence classes of $K$-cycles over
$X$ generated by direct sum and vector bundle modif\/ication
\cite{BD1}. A dif\/ferential $K$-character on  $X$ is a
homomorphism of semi-group $\varphi:\C_*(X) \to \mathbb R/\mathbb
Z$ such that  it is restriction to   the boundary is given  by the
formula
\[
\varphi(\partial(M,E,\phi)) = \int_{M} \phi^\ast (\omega){\rm
Ch}(E){\rm Td}(M)\ \mod \ \mathbb  Z ,
 \]
 where $\omega$ is a closed  form on $X$ with integer $K$-periods \cite{BM}, ${\rm Ch}(E)$ is the Chern form of the connection $\nabla ^E$  and ${\rm Td}(M)$ is the Todd form of the tangent bundle of~$M$. This can be assembled into a group which is denoted by  $\hat K^{*}(X)$ and called the group of dif\/ferential $K$-characters. We showed then that many secondary invariants can be expressed  as a dif\/ferential $K$-characters,  and the group $K^*(X,\mathbb R/\mathbb  Z)$  of $K$-theory of $X$  with coef\/f\/icients in $\mathbb R/\mathbb  Z$ \cite{APS2} is injected in $\hat K^{*}(X)$.

 The aim of this work is to def\/ine the group $\hat K^{*}(\rho)$ of  relative  dif\/ferential $K$-characters associa\-ted with a  smooth map $\rho: Y\to X$ between two  smooth  compact manifolds  $Y$ and $X$ following  \cite{BT,HL} and \cite{HS}. We show that this group f\/its into a short exact sequence as in the non-relative case. The paper is organized as follows:

 In Section~\ref{sec2},
 we def\/ine a group of relative geometric $K$-homology $K_*(\rho)$ adapted to
  this situation and  study some  of its properties. This generalizes the works of Baum--Douglas \cite{BD2} for $Y$  a submanifold of $X$. Section~\ref{sec3} is concerned with the def\/inition and the  study of  the  group  $\hat K^{*}(\rho)$ of  relative  dif\/ferential $K$-characters. An odd relative group $K^{-1}(\rho, \mathbb R/\mathbb Z)$  of $K$-theory with coef\/f\/icients in $\mathbb R/\mathbb Z$ is also def\/ined  here. We  proof  the following short exact sequence
\[
0\to K^{-1}(\rho, \mathbb R/\mathbb Z) \hookrightarrow {\hat
K}^{-1}(\rho) \; {\stackrel{\delta_1}{\rightarrow}}\;
\Omega_{0}^{\rm even}(\rho) \to 0,
\]
where $\Omega_{0}^{\rm even}(\rho)$ is the group of relative
dif\/ferential forms (Def\/inition~\ref{def6}) with integer
$K$-periods. We show then that some secondary geometric invariants
can be expressed in this theory.

 \section[Relative geometric $K$-homology]{Relative geometric $\boldsymbol{K}$-homology}\label{sec2}

 Let $Y$ and $X$ be smooth compact  manifolds and $\rho: Y\to X$ a smooth map. In this section, we def\/ine the relative geometric $K$-homology $K_{*}(\rho)$ for the triple  $(\rho, Y, X)$. This construction generalizes the relative geometric $K$-homology group $K_{*}(X,Y)$ of Baum--Douglas for $Y$ being a closed submanifold of $X$. We recall the def\/inition of the geometric $K$-homology  of a smooth manifold  following the works of Baum and Douglas. This def\/inition  is purely  geometric. For a~complete presentation  see \cite{BD1,BD2} and \cite{RS}.

\begin{definition}\label{def1}
A  $K$-chain over $X$ is a triple $(M,E,\phi)$
 such that:
\begin{itemize}\itemsep=0pt

\item $M$ is a smooth ${\rm Spin}^{c}$-compact  manifold which may
have non-empty boundary  $\partial M $, and with a  f\/ixed
Riemannian structure;

\item $E$ is a Hermitian vector bundle over $M$ with a f\/ixed
Hermitian connection $\nabla ^E$;

\item $\phi : M \to X$ is a smooth map.
\end{itemize}
\end{definition}

Denote that $M$ is not supposed connected and the f\/ibres of $E$
may have dif\/ferent dimensions on the dif\/ferent connected
components of $M$. Two $K$-chains  $(M,E,\phi)$ and
$(M',E',\phi')$  are said to be isomorphic if there exists  a
dif\/feomorphism $\psi:M \to M'$   such that:
\begin{itemize}\itemsep=0pt

\item $\phi'\circ \psi = \phi$;

\item $\psi^*E' \cong E$ as Hermitian bundles over $M$.
\end{itemize}

A $K$-cycle is a $K$-chain $(M,E,\varphi)$ without boundary; that
is  $\partial M = \varnothing$. The boundary
$\partial(M,E,\varphi)$ of the $K$-chain $(M,E,\varphi)$ is the
$K$-cycle $(\partial M, E|_{\partial M}, \varphi|_{\partial M})$.
The set of  $K$-chains is stable under disjoint union.

\subsection[Vector bundle modification]{Vector bundle modif\/ication}

Let $(M,E,\phi)$ be $K$-chain over $X$, and  let $H$ be a ${\rm
Spin}^{c}$-vector bundle over $M$ with even dimensional f\/ibers
and a f\/ixed Hermitian structure. Let $l = M\times \mathbb R$ be
the trivial bundle and $\hat{M} = S(H\oplus l)$  the unit  sphere
bundle. Let $\rho:\hat{M}\to M$  the natural projection. The ${\rm
Spin}^{c}$-structure on $M$ and $H$ induces a ${\rm
Spin}^{c}$-structure on $\hat{M}$.

 Let $\S = \S_{-}\oplus \S_{+}$ be the $\mathbb Z/2\mathbb Z$-grading Clif\/ford module  associated with the ${\rm Spin}^{c}$-structure of $H$. We denote by $H_{0}$ and $H_{1}$ the pullback of $\S_{-}$ and  $\S_{+}$ to $H$. Then $H$ acts on $H_{0}$ and $H_{1}$ by Clif\/ford multiplication: $ H_{0}\stackrel{\sigma}{\rightarrow} H_{1}$.

 The manifold $\hat{M}$ can be thought of as two copies, $B_{0}(H)$ and $B_{1}(H)$, of the unit ball glued together by the identity map of $S(    H)$
\[
\hat{M} = B_{0}(H)\cup_{S(H)}B_{1}(H).
\]

    The vector bundle $\hat{H}$ on $\hat{M}$ is obtained by putting $H_{0}$ on  $B_{0}(H)$ and  $H_{1}$ on $B_{1}(H)$ and then clutching these two vector bundles along  $S(H)$ by the isomorphism $\sigma$.

The  $K$-chain $({\hat M}, {\hat H}\otimes \rho^{*}E, {\hat\phi}=
\rho\circ \phi)$ is called the Bott $K$-chain associated with the
$K$-chain $(M,E,\phi)$ and the ${\rm Spin}^{c}$-vector bundle $H$.

The  boundary of the Bott $K$-chain $({\hat M}, {\hat H}\otimes
\rho^{*}E, {\hat\phi})$ associated with the $K$-chain $(M,E,\phi)$
and the ${\rm Spin}^{c}$-vector bundle $H$ is the Bott $K$-cycle
of the boundary $\partial (M, E, \phi)$ with the restriction of
$H$ to $\partial M$.

\begin{definition}\label{def2}
We denote by $\C_*(X)$ the set of equivalence classes of
isomorphic $K$-cycles over~$X$ up to the following
identif\/ications:

\begin{itemize}\itemsep=0pt

\item we identify the disjoint union  $(M,E,\phi) \amalg
(M,E',\phi)$ with the $K$-cycle $(M,E\oplus E',\phi)$;

\item we identify a $K$-cycle $(M,E,\phi)$ with the Bott $K$-cycle
$({\hat M}, {\hat H}\otimes \rho^{*}E, {\hat\phi})$ associated
with any Hermitian vector bundle $H$ over $M$.
\end{itemize}
\end{definition}

We can easily show that disjoint union then respects these
identif\/ications and makes $\C_*(X)$ into an Abelian semi-group
which splits into $\C_0(X) \oplus \C_1(X)$ with respect to the
parity of the connected components of the manifolds in (the
equivalence classes of) the $K$-cycles.

\begin{definition}\label{def3}
Two $K$-cycles $(M,E,\phi)$ and $(M',E',\phi')$ are bordant if
there exists a $K$-chain $({\overline N}, \E, \psi)$  such that
\[
    \partial({\overline N},\E,\psi)   \quad \mbox{is
isomorphic   to} \ \
    (M,E,\phi) \amalg (-M', E',\phi'),
\]
    where $-M'$ is $M'$ with the  ${\rm Spin}^{c}$-structure reversed~\cite{BD1} .
\end{definition}

The  above bordism relation induces a well def\/ined equivalence
relation on $\C_*(X)$ that we denote by $\sim_{\partial}$. The
quotient $\C_*(X)/\sim_{\partial}$ turns out to be an Abelian
group for the disjoint union. The inverse of $(M, E,\phi)$ is
$(-M, E,\phi)$.

\begin{definition}[Baum--Douglas]\label{def4}
The quotient group of $\C_*(X)$ by the equivalence relation
$\sim_{\partial}$ is denoted by $K_{*}(X)$ and is called the
geometric $K$-homology group of $X$. It can be decomposed into
\[
    K_{*}(X) = K_{0}(X) \oplus K_{1}(X).
\]
\end{definition}

 A smooth map $\varphi:Y \to X$
induces a group morphism
\[
    \varphi_*:  K_{*}(Y) \to  K_{*} (X),
\]
given by $\varphi_*(f)(M,E,\phi) = f(M,E,\varphi \circ \phi)$. The
$ K_{*}$ is a covariant functor from the category of smooth
compact manifolds and smooth maps to that of Abelian groups and
group homomorphisms.

In the same way we can form a semi-group $\L_*(X)$ out of
$K$-chains $(\overline N,\E,\psi)$, say with the same def\/inition
as $\C_*(X)$ and the \textit{boundary}
\[
\partial({\overline N},\E,\psi) = (\partial \overline{N}, \E|_{\partial \overline{N}},\psi \circ i),
\]
 \noindent where $i:\partial \overline {N}\hookrightarrow {\overline N}$. This gives  a well def\/ined  map
\[
    \partial: \L_*(X) \to \C_*(X) \subset \L_*(X).
 \]

The Hermitian structure of the complex vector bundle
$\E|_{\partial \overline{ N}}$ is the restricted one.

The group of $K$-cochains with coef\/f\/icients in $\mathbb Z$ denoted by $\L^*(X)$ is the group of semi-group homomorphisms $f$ from $\L_*(X)$ to $\mathbb Z$.
On the group $\L^*(X)$ there is  a coboundary map
def\/ined by transposition
\[
    \delta(f)({\overline N}, \E, \psi) = f(\partial (\overline N, \E, \psi)).
\]

   The set of $K$-cocycles is the subset $\C^*(X)$ of $\L^*(X)$ of those $K$-cochains that
vanish on boundaries, i.e. the kernel of $\delta$. The set of
$K$-coboundaries   is the image of $\delta$ in $\L^*(X)$.

\subsection[The relative geometric group $K_{*} (\rho)$]{The relative geometric group $\boldsymbol{K_{*} (\rho)}$}\label{sec2.1}

Let $Y$ and $X$ be  smooth compact manifolds and $\rho: Y\to X$ a
smooth map.

 The set $\L_{*}(\rho)$ of relative $K$-chains associated with the triple $(\rho,Y,X)$ is by def\/inition
\[
\L_{*+1}(\rho) = \L_{*+1}(X)\times \L_{*}(Y).
\]
The boundary $\partial:\L_{*+1}(\rho)\to \L_{*}(\rho)$ is given by
\[
\partial(\sigma,\tau) = (\partial\sigma + \rho_{*}\tau, -\partial\tau).
\]

We will denote by $\C_{*}(\rho)$ the set of relative $K$-cycles in
$\L_{*}(\rho)$, i.e., the kernel of $\partial$. A  $K$-cycle in
$\L_{*}(\rho)$ is then a pair $(\sigma,\tau)$ where $\tau$ is a
$K$-cycle over $Y$ and $\sigma$ is $K$-chain over $X$ with
boundary in the image of $\rho_{*}:\C_{*}(Y)\to \C_{*}(X)$.  The
set $\C_{*}(\rho)$ is a semi-group for the sum
\[
(\sigma,\tau)+ (\sigma',\tau') =  (\sigma \amalg\sigma',\tau
\amalg \tau'),
\]
where $\amalg$ is the disjoint union. We say that two relatives
$K$-cycles$(\sigma,\tau)$ and $(\sigma',\tau')$  are bordant and
we write $(\sigma,\tau)\sim_{\partial} (\sigma',\tau')$   if there
exists a relative $K$-chain $(\overline{\sigma},\overline{\tau})$
such that
\[
\partial (\overline{\sigma},\overline{\tau}) = (\sigma,\tau)+ (-\sigma',-\tau'),
\]
where $-x$ denotes the relative $K$-cycle $x$ with the reversed
${\rm Spin}^{c}$-structure of the underlying manifold.

\begin{definition}\label{def5}
 The relative geometric $K$-homology group  denoted by $ K_{*}(\rho)$ is   the quotient group $\C_{*}(\rho)/\sim_{\partial}$.
 \end{definition}

 The inverse of the $K$-cycle $x$ is $-x$. The equivalence relation on the relative $K$-cycle $(\sigma, \tau)$ preserves the dimension modulo 2 of the $K$-cycles  $\sigma$ and $\tau$. Hence, there is a direct sum decomposition
\[
 K_{*}(\rho) =  K_{0}(\rho)\oplus  K_{1}(\rho).
 \]

 The construction of the group $ K_{*}(\rho) $ is functorial in the sense that
 for a commutative diagram
$$ \begin{array}{clrr} %
       Y &  \stackrel{\rho}{\longrightarrow}& X &  \\
       \Big\downarrow{f}{}& &\Big\downarrow{}{g}\\
          Y'  &   \stackrel{\rho'}{\longrightarrow}& X' \cr
 \end{array}
$$
the map $F_{*}=(f_*,g_*):\L_{*}(\rho)\to \L_{*}(\rho')$ is
compatible with the equivalence relation on the relative
$K$-cycles and induces a homomorphism from $K_{*}(\rho)$ to
$K_{*}(\rho')$. As in the homology theory, we have the  long exact
sequence for the triple $(\rho,Y,X)$
$$
\begin{array}{clrrr}
{K}_{0} (Y) &\stackrel{ \rho_{*}}{\longrightarrow}& {K}_{0}(X)&\stackrel{\varsigma_{*}}{\longrightarrow}&  {K}_{0}(\rho) \\
\Big\uparrow{\partial}{}&  &   &  &\Big\downarrow{}{\partial}\\
{K}_{1}(\rho)&\stackrel{\varsigma _{*}}{\longleftarrow}&
{K}_{1}(X)&\stackrel{\rho_{*}}{\longleftarrow}&{K}_{1} (Y)\cr
 \end{array}
$$

The boundary map  $\partial$ associates to  a relative $K$-cycle
$(\sigma,\tau)$ the cycle $\tau$ whose image  $\rho_{*}\tau$  is
a~boundary in $X$ and $\varsigma_{*}(\sigma) = (\sigma,0)$. The
exactness of the diagram is an easy check.

There is a  dif\/ferential $\delta$ on the group
$\L^{*}(\rho) = {\rm Hom}( \L_{*}(\rho),\mathbb Z)$
 given by
\[
\delta(h,e) = (\delta h, \rho^{*} h - \delta e).
\]

The relative Baum--Douglas $K$-group is
\[
K^{*}(\rho)= \frac{{\rm ker}(\delta :\L^{*}(\rho)\rightarrow \L^{*
+ 1}(\rho))}{{\rm Im}( \delta :\L^{* - 1}(\rho) \rightarrow
\L^{*}(\rho))}.
\]

\begin{remark} 
The relative topological $K$-homology group $K_{*}^{t}(\rho)$ can be constructed
in the same way for normal topological spaces~$X$ and~$Y$, and $\rho:Y\to X$ is a continuous map.   Let $ K_{*}^{t}(X,Y)$ be the relative topological  $K$-homology group def\/ined by Baum--Douglas in \cite{BD2} for $Y\subset X$ is a~closed subset of a $X$. We can easily show that $ K_{*}^{t}(X,Y) = K_{*}^{t}(\rho)$, where   $\rho$ is the inclusion of~$Y$ in~$X$.
\end{remark}

\section[Relative differential $K$-characters]{Relative dif\/ferential $\boldsymbol{K}$-characters}\label{sec3}

This section is concerned with the def\/inition and the  study the
notion of relative dif\/ferential $K$-characters \cite{BM}. This
is a  $K$-theoretical version  of the works of \cite{BT,HL} and
\cite{HS}.

Let $X$ be  a smooth compact manifold. The graded dif\/ferential
complex of real dif\/ferential forms on the manifold $X$ will be
denoted by
\[
\Omega^*(X)=\oplus_{k\geq 0} \Omega^k(X), \qquad \Omega^k(X)
\stackrel{d}{\rightarrow} \Omega^{k+1}(X)\qquad
 \mbox{with} \quad d^2=0,
\]
where $d$ denotes the de Rham dif\/ferential on $X$.

Furthermore, we denote  by  $\Omega_{0}^*(X)$ the subgroup of
closed forms on the manifold $X$ with integer
$K$-periods~\cite{BM}.

In the remainder of this section we f\/ix $\rho: Y\to X$ a smooth
map and we consider the complex
\[
\Omega^*(\rho) = \Omega^*(X) \times \Omega^{*-1}(Y)
\]
with dif\/ferential $\delta(\omega,\theta) = (d \omega,
\rho^{*}\omega - d\theta).$

We can view $\Omega^*(\rho)$ as a subgroup of the the group ${\rm
Hom}(\L_*(\rho),\mathbb R)$ via integration
\[
(\omega,\theta)(\sigma,\tau) = \omega(\sigma) + \theta(\tau),
\]
where  for $\sigma = (M,E,f)$ and $\tau = (N,F,g)$
\[
\omega(\sigma) = \int_{M}f^{*}(\omega){\rm Ch}(E){\rm Td}(M)
\qquad \mbox{and} \qquad \theta(\tau) = \int_{N}g^{*}(\theta){\rm
Ch}(F){\rm Td}(N).
\]

 Let
\[
    j:\Omega^*(\rho) \to {\rm Hom}(\L_*(\rho),\mathbb R)
\]
such that
\[
 j(\omega,\theta)(\sigma,\tau) = \omega(\sigma) + \theta(\tau).
 \]

\begin{definition}\label{def6}
Let $(\omega,\theta)\in \Omega^*(\rho)$ be a pair of real
dif\/ferential forms.
\begin{enumerate}\itemsep=0pt
\item[(i)] The set of $K$-periods of $(\omega,\theta)$ is the
subset of $\mathbb R$ image of the map $j(\omega,\theta)$
restricted to~$\C_*(\rho)$.

\item[(ii)] We denoted by $\Omega^*_0(\rho)$ the set of
dif\/ferential forms $(\omega,\theta)$  of integer   $K$-periods.
    \end{enumerate}
\end{definition}

 $\Omega^*_0(\rho)$ is an Abelian group for the sum of dif\/ferential forms.

\begin{lemma}\label{lemma1}
 Let $(\omega,\theta)\in \Omega^*_0(\rho)$. Then

\begin{enumerate}\itemsep=0pt
\item[\rm 1)] $\delta(\omega,\theta) = 0 $ in the complex
$\Omega^{*}(\rho)$;

\item[\rm 2)] $\omega\in\Omega^{*}_{0}(X)$.
\end{enumerate}
\end{lemma}

\begin{proof} 1) For $(\omega,\theta)\in \Omega^{*}_{0}(\rho)$ and $ \tau =  (N,F,g)\in\L_{*-1}(Y)$, we have
\begin{gather*}
\rho^{*}\omega(\tau) - d \theta(\tau)=\int_{N}(g^{*}(\rho^{*}(\omega)) - g^{*}(d \theta)){\rm Ch}( F) {\rm Td}(N) \\
\phantom{\rho^{*}\omega(\tau) - d \theta(\tau)}{}
 =   \int_{N}g^{*}(\rho^{*}(\omega)){\rm Ch}( F) {\rm Td}(N) - \int_{\partial N} g^{*}( \theta){\rm Ch}( F) {\rm Td}(N) \\
\phantom{\rho^{*}\omega(\tau) - d \theta(\tau)}{}
 =   \int_{N}(\rho\circ g)^{*}(\omega){\rm Ch}( F) {\rm Td}(N) - \int_{\partial N} g^{*}( \theta){\rm Ch}( F) {\rm Td}(N)\\
\phantom{\rho^{*}\omega(\tau) - d \theta(\tau)}{} =
(\omega,\theta)(\rho_{*}\tau,-\partial \tau).
\end{gather*}

Since  $(\omega,\theta)\in \Omega^{*}_{0}(\rho)$ and
$(\rho_{*}\tau,-\partial \tau) = \partial(0,\tau) $ is a relative
$K$-cycle, the value $(\omega,\theta)(\rho_{*}\tau,-\partial
\tau)$ is entire. Lemma 3 of \cite{BM} implies that
$\rho^{*}\omega - d \theta = 0$. On the other hand, for any
$K$-chain $\sigma \in \L(X)$, we have $d\omega (\sigma ) =
(\omega,\theta)(\partial\sigma,0)$. Since $(\partial\sigma,0) $ is
a relative $K$-cycle, it follows for the same raison that
$d\omega = 0$.

2) Let $\sigma = (M,E,f)\in\C_{*}(X)$. We have
\[
\int_{M} f^\ast (\omega){\rm Ch}(E){\rm Td}(M) = (\omega,\theta)(\sigma,0).
\]

Since $(\omega,\theta)$ has integer $K$-periods and $(\sigma,0)$
is a relative $K$-cycle, the right hand-side is entire.~~
\end{proof}

\begin{example}  Any pair $(\omega,\theta)\in \Omega^*(\rho)$ of exact dif\/ferential forms is obviously in  $\Omega_0^*(\rho)$.
\end{example}

\begin{remark}\label{rem2}
 We can easily deduce from the proof of the previous lemma that an element $(\omega,\theta)\in\Omega^{*}(\rho)$ with  entire values on all $K$-chains is necessarily trivial.
\end{remark}

\begin{definition}\label{def7} \ \null
\begin{enumerate}\itemsep=0pt
\item[(i)] A relative dif\/ferential $K$-character for the smooth
map $\rho:Y\to X$  is a homomorphism of semi-group
\[
    f: \C_*(\rho) \to \mathbb R/\mathbb Z
\]
such that $f(\partial(\sigma,\tau)) =
[(\omega,\theta)(\sigma,\tau)] $ for some
$(\omega,\tau)\in\Omega_{0}^*(\rho)$ and for all relative
$K$-chain $(\sigma,\tau)\in\L_{\ast}(\rho)$, where $[\alpha]$
denote the class in $\mathbb R/\mathbb Z$ of the number $\alpha$.

\item[(ii)] The set of relative dif\/ferential $K$-characters is
denoted by ${\hat K}^{*}(\rho)$. It is naturally $\mathbb
Z/2\mathbb Z$-graded
\[
{\hat K}^{*}(\rho) = {\hat K}^{0}(\rho) \oplus {\hat K}^{1}(\rho).
\]
\end{enumerate}
\end{definition}

 Let $f$ be a relative  dif\/ferential $K$-character for the smooth map $\rho: Y\to X$. We  deduce  from Remark~\ref{rem2} that the pair  of forms $(\omega,\theta)$ associated to $f$ in Def\/inition~\ref{def7} is unique. It will be denoted  by  $\delta_1(f)$.
We thus have a group morphism
\[
    \delta_1: {\hat K}^{*}(\rho) \to \Omega_{0}^*(\rho),
\] which is odd for the grading.

\begin{example} An interesting situation  is obtained by dif\/ferential forms.
If $(\omega,\theta)\in \Omega^*(X)\times \Omega^{*-1}(Y)$ is any
pair of real dif\/ferential forms, then we def\/ine
$f_{(\omega,\theta)}$ by letting for $\sigma = (M,E,f)$ and $\tau
= (N,F,g)$
\[
    f_{(\omega,\theta)}(\sigma,\tau) = \left[ \int_{M}f^{*}(\omega){\rm Ch}(E){\rm Td}(M)\right]   + \left[\int_{N}g^{*}(\theta){\rm Ch}(F){\rm Td}(N)\right].
    \]

We have
\[
\delta_{1}(f_{(\omega,\theta)}) = (d\omega, \rho^{*}\omega -
d\theta).
\]
\end{example}

\begin{example} Suppose  $Y$ be submanifold of $X$ and $\rho: Y \hookrightarrow X$ is  a smooth inclusion. Let $\omega\in \Omega^*(X)$ with trivial restriction to $Y$ and $\bar{ f_{\omega}} \in \hat{K}(X)$ -- the associated  dif\/ferential $K$-character~\cite{BM}. Let $\psi\in \hat{K}(Y)$  be any  dif\/ferential $K$-character on $Y$. We have  $\bar{ f_{\omega}}(\L_{*}Y) = 0$. The map  $\phi_{\omega,\psi}$ def\/ined on $\C_{*}(\rho)$ by
\[
\phi_{\omega,\psi}(\sigma,\tau) = \bar{ f_{\omega}} (\sigma) +
\psi)(\tau)
 \]
 is a relative dif\/ferential $K$-character with $\delta_{1}(\phi_{\varphi,\psi}) = (d\omega, - \delta_{1}(\psi)).$
\end{example}

 \subsection[Relative $\mathbb R/\mathbb Z$-$K$-theory]{Relative $\boldsymbol{\mathbb R/\mathbb Z$-$K}$-theory}\label{sec3.1}

Let $X$ be  a smooth  manifold,
 $E$  a Hermitian vector bundle on $X$ and $\nabla^{E}$  a Hermitian connection on $E$. The geometric Chern form ${\rm Ch}(\nabla^{E})$  of $\nabla^{E}$  is the closed  real  even dif\/ferential form given~by
\begin{gather*}
{\rm Ch}(\nabla^{E}) = {\rm
tr}\,e^{-\frac{(\nabla^{E})^{2}}{2i\pi}}.
\end{gather*}

The cohomology class of ${\rm Ch}(\nabla^{E})$ does not depend on
the choice of the connection~$\nabla^{E}$~\cite{Lott}.
Let~$\nabla_{1}^{E}$ and  $\nabla_{2}^{E}$ be two Hermitian
connections on $E$. There is a well def\/ined  Chern--Simons
form~\cite{Lott} ${\rm CS}(\nabla_{1}^{E},\nabla_{2}^{E}) \in
\frac{\Omega^{\rm odd}(X)\otimes \mathbb C}{{\rm Im}(d)}$  such
that
\[
d{\rm CS}(\nabla_{1}^{E},\nabla_{2}^{E}) = {\rm
Ch}(\nabla_{1}^{E}) - {\rm Ch}(\nabla_{2}^{E}),
\]
and
\[
{\rm  CS}(\nabla_{1}^{E},\nabla_{3}^{E}) = {\rm
CS}(\nabla_{1}^{E},\nabla_{2}^{E}) + {\rm CS}
(\nabla_{2}^{E},\nabla_{3}^{E}).
\]

Given a short exact sequence of complex Hermitian vector bundles
on $X$
\[
0\to  E_{1} \; {\stackrel{i}{\rightarrow }}\;  E_{2} \;
{\stackrel{j}{\rightarrow }}\; E_{3}\to 0,
\]
and choose a splitting map $s:E_{3}\rightarrow E_{2}$. Then
$i\oplus s: E_{1}\oplus E_{3}\rightarrow E_{2}$ is an isomorphism.
For~$\nabla^{E_{1}}$, $\nabla^{E_{2}}$ and $\nabla^{E_{3}}$ are
Hermitian connection on $E_{1}$, $E_{2}$ and $ E_{3}$
respectively, we set
\[
{\rm CS}(\nabla^{E_{1}}, \nabla^{E_{2}},\nabla^{ E_{3}}) = {\rm
CS}( (i\oplus s)^{*}\nabla^{E_{2}}, \nabla^{E_{1}}\oplus\nabla^{
E_{3}}).
\]
 The form  ${\rm CS}(\nabla^{E_{1}}, \nabla^{E_{2}},\nabla^{ E_{3}})$ is independent of the choice of the splitting map $s$ and
\[
d{\rm CS}(\nabla^{E_{1}}, \nabla^{E_{2}},\nabla^{ E_{3}}) = {\rm
Ch}(\nabla_{2}^{E}) - {\rm Ch}(\nabla_{1}^{E}) - {\rm
Ch}(\nabla_{3}^{E}).
\]

\begin{definition}\label{def8}
 Let  $X$  be a smooth manifold. A  $\mathbb R/\mathbb Z$-$K$-generator of $ X$ is a quadruple
\[
\E = (E,h^{E},\nabla^{E},\omega),
\]
where
 $E$ is a complex vector bundle on $X$,
 $h^{E}$ is a positive def\/inite Hermitian metric on $E$,
 $\nabla^{E}$~is a Hermitian connection on $E$,
 $\omega\in \frac{\Omega^{\rm odd}(X)}{{\rm Im}(d)}$ which  satisf\/ies  $d\omega = {\rm Ch}(\nabla^{E}) - {\rm rank}(E)$, where ${\rm rank}(E)$ is the rank of $E$.
\end{definition}

  An $\mathbb R/\mathbb Z$-$K$-relation is given by three $\mathbb R/\mathbb Z$-$K$-generators $\E_{1}$,  $\E_{2}$,  $\E_{3}$, along with a short exact sequence of Hermitian vector bundles
\[
0\to  E_{1}\; {\stackrel{i}{\rightarrow }}\; E_{2} \;
{\stackrel{j}{\rightarrow}}\;
 E_{3}\to 0,
\]
such that $\omega_{2} = \omega_{1} + \omega_{3} + {\rm
CS}(\nabla^{E_{1}},\nabla^{E_{2}},\nabla^{E_{3}})$.

\begin{definition}[\cite{Lott}]\label{def9}
 We denote by $MK(X,\mathbb R/\mathbb Z)$ the quotient of the free group ge\-ne\-ra\-ted by the $\mathbb R/\mathbb Z$-$K$-ge\-ne\-ra\-tors and $\mathbb R/\mathbb Z$-$K$-relation   $\E_{2} = \E_{1} + \E_{3}$. The group $K^{-1}(X,\mathbb R/\mathbb Z)$ is the subgroup of $MK(X,\mathbb R/\mathbb Z)$ consisting of elements of virtual rank zero.
 \end{definition}

 The elements of $K^{-1}(X,\mathbb R/\mathbb Z)$ can be described by $\mathbb Z/2\mathbb Z$-graded cocycles \cite{Lott}, meaning quad\-rup\-les $ (E_{\pm},h^{E_{\pm}},\nabla^{E_{\pm}},\omega)$ where
 $E = E_{+}\oplus E_{-}$   is a  $\mathbb Z/2\mathbb Z$-graded complex vector bundle on $X$,
 $h^{E} = h^{E_{+}}\oplus h^{E_{-}}$ is a Hermitian metric on $E$,
 $\nabla^{E}=\nabla^{E_{+}}\oplus \nabla^{E_{-}}$ is a Hermitian connection on~$E$,
$\omega\in \frac{\Omega^{\rm odd}(X)}{{\rm Im}(d)}$ and
satisf\/ies $d\omega = {\rm Ch}(\nabla^{E})= {\rm
Ch}(\nabla^{E_{+}}) - {\rm Ch}(\nabla^{E_{-}})$.

    We consider now  two smooth compact  manifolds $Y$ and $X$. Let $\rho: Y\to X$ be a smooth  map and let the exact sequence
$$
 \begin{array}{ccccc}
{K}^{0} (Y,\mathbb R/\mathbb Z) &\stackrel{ \rho_{0}^{*}}{\longleftarrow}& {K}^{0}(X,\mathbb R/\mathbb Z)&\stackrel{\varsigma_{0}^{*}}{\longleftarrow}  & {\rm Hom}({K}_{0}(\rho),\mathbb R/\mathbb Z)  \\
 \Big\downarrow{\partial_{0}^{*}}{}&  &  & &\Big\uparrow{}{\partial_{1}^{*}}\\
{\rm Hom}({K}_{1}(\rho),\mathbb R/\mathbb Z)
&\stackrel{\varsigma_{1}^{*}}{\longrightarrow}& {K}^{-1}(X,\mathbb
R/\mathbb Z)&\stackrel{\rho_{1}^{*}}{\longrightarrow} &
{K}^{-1}(Y,\mathbb R/\mathbb Z)
 \end{array}
$$
obtained from the one in p.~4 and after identif\/ication of the
groups $K^{*} (Y,\mathbb R/\mathbb Z)$ and ${\rm Hom}(K_{*} (Y)$,
$\mathbb R/\mathbb Z)$ following Proposition~4 of~\cite{BM}.

 We denote by   ${\bar L} ^{-1}(\rho, \mathbb R/\mathbb Z)$ the subgroup of ${\rm Hom}({K}_{1}(\rho),\mathbb R/\mathbb Z)$  image of the morphism $\partial_{0}^{*}$ and ${\bar K} ^{-1}(\rho, \mathbb R/\mathbb Z)$ the subgroup of ${K}^{-1}(X,\mathbb R/\mathbb Z)$ the kernel of the morphism $\rho_{1}^{*}$.

\begin{definition}\label{def10}  Let $Y$ and $X$ be two smooth compact  manifolds and $\rho: Y\to X$ a smooth map.  The group   $K^{-1}(\rho, \mathbb R/\mathbb Z)$  is by def\/inition the product of the groups ${\bar L^{-1}}(\rho, \mathbb R/\mathbb Z)$  and ${\bar K} ^{-1}(\rho, \mathbb R/\mathbb Z)$
\[
K^{-1}(\rho, \mathbb R/\mathbb Z) = {\bar L^{-1}} (\rho, \mathbb
R/\mathbb Z)\times {\bar K} ^{-1}(\rho, \mathbb R/\mathbb Z).
\]
\end{definition}

 \begin{proposition}\label{prop1} The groups  $K^{-1}(\rho, \mathbb R/\mathbb Z)$ and  ${\rm Hom}(K_{-1}(\rho), \mathbb R/\mathbb Z)$ are isomorphic.
 \end{proposition}

 \begin{proof} Since the image of $\varsigma_{1}^{*}$ is the kernel of $\rho_{1}^{*}$, it is enough to show that the short exact sequence
\[
 0\rightarrow {\bar L^{-1}}(\rho, \mathbb R/\mathbb Z)\hookrightarrow  {\rm Hom}({K}_{1}(\rho),\mathbb R/\mathbb Z)\; {\stackrel{\varsigma _{1}^{*}}{\to}} \; {\bar K}^{-1}(\rho,\mathbb R/\mathbb Z)\rightarrow 0
 \]
 is split. Let  $\E$ be an element of ${\bar K}^{-1}(\rho, \mathbb R/\mathbb Z)$ and $(E_{\pm},h^{E_{\pm}},\nabla^{E_{\pm}},\omega)$ be a relative $\mathbb Z/2\mathbb Z$-graded cocycle associated to $\E$. Let $ (\sigma,\tau)$  be a relative  $K$-cycle in $\C_{-1}(\rho)$.  For $\sigma =(M, E, \phi)$ and $\tau =(N,F,\psi)$ we set
\[
 \alpha(\E)(\sigma,\tau)) = \bar{\eta}_{\phi^{*}E_{+}\otimes E} -\bar{\eta}_{\phi^{*}E_{-}\otimes E}
 -\bar{f}_{\omega}(\sigma),
 \]
 where the notation $\bar{\eta}_{V} = \frac{\eta(D_V) + \dim \ker D_V}{2}$ $(\!\!\!\!\mod  \mathbb Z$) is the reduced eta invariant~\cite{APS2,APS3} of Atiyah--Patodi--Singer  of the Dirac operator associated to the ${\rm Spin}^c$-structure of $M$ with coef\/f\/icients in the vector bundle $V$ \cite{APS1} and
\[
\bar{f}_{\omega}(M,E,\phi) =\left[\int_{M}\phi^{* }(\omega) {\rm
Ch}(E){\rm Td}(M)\right].
\]
   Let us check that $\alpha(\E)(\partial (\sigma,\tau)) = 0$ in $\mathbb R/\mathbb Z$ for any $K$-chain $\sigma$ over $X$ and any $K$-chain $\tau$ over $Y$. Recall that  $\partial (\sigma,\tau)= (\partial \sigma + \rho^{*}\tau, -\partial\tau)$. Furthermore, the invariant $\bar{\eta}$ and $\bar{f}_{\omega}$ def\/ines  $K$-cochains over $X$ \cite{BM}. We have then
\[
\alpha(\E)(\partial (\sigma,\tau))= \alpha(\E)(\partial\sigma,
-\partial\tau) + \alpha(\E)(\rho^{*}\tau,0).
\]

The index theorem of APS (see \cite{APS1,APS2,APS3}) implies that
\[
\bar{\eta}_{(\phi^{*}E_{+}\otimes E)|\partial M}
-\bar{\eta}_{(\phi^{*}E_{-}\otimes E)|\partial M}
 -\bar{f}_{d\omega}(\sigma) = {\rm ind}(D_{+}\otimes \phi^{*}E_{+} \otimes E) - {\rm ind}(D_{+}\otimes \phi^{*}E_{-} \otimes E) ,
 \]
 is entire, where ${\rm ind}(D_{+}\otimes \phi^{*}E_{\pm} \otimes E)$ is the index of the Dirac type operator associated to the ${\rm Spin}^{c}$-structure of $M$ with coef\/f\/icients in the bundle $ \phi^{*}E_{\pm} \otimes E$. On the  other hand, we have
\[
\alpha(\E)(\rho^{*}\tau,0)) = \alpha(\rho^{*}\E)(0,\tau)=  0.
\]

This construction def\/ines a homomorphism
$\alpha:\bar{K}^{-1}(\rho, \mathbb R/\mathbb Z)\to {\rm
Hom}(K_{-1}(\rho), \mathbb R/\mathbb Z)$ which is a split of
$\varsigma _{1}^{*}$. In fact, let us   consider the following
commutative diagram
$$
  \begin{array}{ccc}
       \bar{K}^{-1}(\rho, \mathbb R/\mathbb Z)  & \stackrel{\alpha}{\longrightarrow} & {\rm Hom}(K_{-1}(\rho), \mathbb R/\mathbb Z)\\
        \Big\downarrow{i^{*}}{} && \Big\downarrow{\varsigma ^{*}_1}{}\\
        {K}^{-1}(X,\mathbb R/\mathbb Z)   & \stackrel{\alpha _{X}}{\longrightarrow} & {\rm Hom}(K_{-1}(X), \mathbb R/\mathbb Z)
  \end{array}
$$
where $i^{*}$ is the embedding of $\bar{K}^{-1}(\rho, \mathbb
R/\mathbb Z)$ in $K^{-1}(X,\mathbb R/\mathbb Z)$ and  $\alpha
_{X}$ is the restriction of $\alpha $ to~$\C_{*}(X)\times \{0\}$
which is  an isomorphism \cite{Lott}. For any  $\E\in
\bar{K}^{-1}(\rho, \mathbb R/\mathbb Z)$, we have $\varsigma
^{*}_1(\alpha(\E)) = i^{*}(\E) = \E $.
\end{proof}

\begin{theorem}\label{theorem1}
The following sequence is exact:
\[
0\to K^{-1}(\rho, \mathbb R/\mathbb Z) \hookrightarrow {\hat
K}^{-1}(\rho) \; {\stackrel{\delta_1}{\to}}\; \Omega_{0}^{\rm
even}(\rho) \to 0.
\]
\end{theorem}

\begin{proof} From Proposition \ref{prop1},  $K^{-1}(\rho, \mathbb R/\mathbb Z)$ is isomorphic to
${\rm Hom}( K_{-1}(\rho),\mathbb R/\mathbb Z )$ which obviously
injects in ${\hat K}^{-1}(\rho)$ with trivial $\delta_1$. It is
clear that a relative dif\/ferential $K$-character $f$ with
$\delta_1(f) = 0$, induces a homomorphism from $K_{-1}(\rho)$ to
$\mathbb R/\mathbb Z$. Hence the sequence is exact at
$K^{-1}(\rho)$. It remains to show the surjectivity of $\delta_1$.

Let $(\omega,\theta)\in \Omega_{0}^{\rm even}(\rho)$ and
$f_{\omega,\theta}: \L_{*}(\rho)\rightarrow \mathbb R/\mathbb Z $
def\/ined by
\[
f_{\omega,\theta}(\sigma,\tau) = \overline{f_{\omega}(\sigma)}+
\overline{f_{\theta}(\tau)}.
\]

 The map $f_{\omega,\theta}$ is trivial on $\C_{-1}(\rho)$. Therefore, we def\/ine an element $\chi\in {\rm Hom}(\B_{-1}(\rho),\mathbb R/\mathbb Z )$  by setting
\[
\chi(\partial (\sigma,\tau)) = f_{\omega,\theta}(\sigma,\tau),
\]
where $\B_{-1}(\rho)$ is the image of the boundary map
$\partial:\L_{0}(\rho)\to \C_{-1}(\rho)$.

 Since $\mathbb R/\mathbb Z$ is divisible, $\chi$ can be extended  to a relative dif\/ferential $K$-character ${\overline\chi}: \C_{-1}(\rho) \to \mathbb R/\mathbb Z$ with $\delta_{1}( {\overline\chi}) = (\omega,\theta).$
 \end{proof}

\subsection{Application} Let $G$ be an almost connected Lie group and $M$ be a smooth compact manifold. Let $\pi: Y\to M$ be a compact principal $G$-bundle with connection $\nabla$. We denote by $I^{*}(G)$  the ring of  symmetric multilinear real functions on the Lie algebra of $G$ which are invariant under the adjoint action of~$G$~\cite{ChS}. Let $\Omega$ be the curvature of $\nabla$. For any $P\in I^{*}(G)$, there is a well def\/ined closed form~$P(\Omega)$ on $M$. The pullback $\pi^{*}P(\Omega)$ is an exact form on $Y$ \cite {ChS}. For  $P\in I^{*}(G)$, let~$TP(\nabla)$ be such that $\pi^{*}P(\Omega) = dTP(\nabla)$. The form $\omega = (\pi^{*}P(\Omega), dTP(\nabla))$ is a closed form in the complex $(\Omega^{*}(\pi),\delta)$. The relative dif\/ferential $K$-character $f_{\omega}$ has a trivial $\delta_{1}$ and def\/ines consequently an element of the  group  $K^{-1}(\pi, \mathbb R/\mathbb Z)$. This gives an additive  map from $I^{*}(G)$ to $K^{-1}(\pi, \mathbb R/\mathbb Z)$  which can be looked as a home of secondary geometric invariants of  the  principal $G$-bundle with connection $(M,Y,\nabla)$ analogous to the Chern--Weil theory.

\subsection*{Acknowledgements}
I should like to thank the referees for their very helpful suggestions and important remarks.

\pdfbookmark[1]{References}{ref}

\LastPageEnding

\end{document}